\documentclass[12pt]{amsart}

\usepackage{charter}
\usepackage{setspace}
\usepackage{macros}
\standardsettings

\renewcommand{\df}{:=}
\newcommand{\cardQ}{r}
\newcommand{\poly}{p}
\renewcommand{\path}{P}

\newcommand{\mult}{\times}

\newcommand{\bothdim}{{\pdim,\qdim}}

\renewcommand{\prec}{\preceq}
\renewcommand{\succ}{\succeq}

\draftfalse

\begin{document}
\title[The matrix version of Baker's conjecture]{A proof of the matrix version of Baker's conjecture\\ in Diophantine approximation}

\authortushar\authordavid


\begin{abstract}
We prove that the matrix analogue of the Veronese curve is strongly extremal in the sense of Diophantine approximation, thereby resolving a question posed by Beresnevich, Kleinbock, and Margulis ('15) in the affirmative.
%
%
\end{abstract}
\maketitle

\section{Introduction}
The origins of metric Diophantine approximation on manifolds (also known as Diophantine approximation with dependent quantities) may be traced back to Kurt Mahler's profound conjecture in transcendence theory, \cite{Mahler2}, that may be translated thus: Lebesgue almost every point on the Veronese curve $\{(x,x^2,\ldots,x^n) : x\in\R\}$ is not very well approximable by rational vectors. Mahler's conjecture was resolved three decades later by Vladimir Sprind\v zuk \cite{Sprindzuk2}, who went on to conjecture \cite{Sprindzuk3} that any analytic nondegenerate submanifold of $\R^n$ is strongly extremal, viz. that Lebesgue almost every point is not very well multiplicatively approximable by rational vectors. A special case of Sprind\v zuk's conjecture was conjectured earlier by Alan Baker, \cite[p.96]{Baker_book}, viz. that the Veronese curve is strongly extremal. We refer the reader to \cite{BernikDodson} for more regarding this history and allied developments.

After a slew of partial results by a number of authors, Sprind\v zuk's conjecture was finally resolved by Dmitry Kleinbock and Grigory Margulis in their landmark 1998 Annals paper \cite{KleinbockMargulis2}, where they translated the problem into dynamical terms and successfully leveraged quantitative non-divergence estimates for unipotent flows on homogenous spaces. There has since been an intense study of Diophantine extremality, both in the case of simultaneous approximation \cite{KLW, StratmannUrbanski1, Urbanski} and in the matrix approximation framework \cite{Kleinbock6, KMW, BKM, ABRdS, DFSU_GE1}.
In particular, in a joint work with Junbo Wang, Kleinbock and Margulis generalized their result to the matrix case by proving that any strongly nonplanar manifold in the space of matrices is strongly extremal \cite[Theorem 2.1]{KMW} (cf. Definition \ref{definitionstronglynonplanar} below). 
Strong nonplanarity is not only a stronger requirement than nondegeneracy, but also one that is more difficult to verify. On the other hand, it is necessary since Kleinbock--Margulis--Wang produced examples of nondegenerate manifolds in the space of matrices which are not extremal \cite[Conclusion 2 on p.22]{KMW}.

It is natural to ask whether the following matrix version of Mahler's (resp. Baker's) conjecture holds: If $m,n\in\N$ are fixed, then for Lebesgue almost every matrix $X\in \MM_{m,m}$, the matrix $X\oplus X^2\oplus \cdots\oplus X^n\in \MM_{m,mn}$ is not very well approximable (resp. not very well multiplicatively approximable). The matrix analogue of Mahler's conjecture was proven by Kleinbock \cite[\63.3]{Kleinbock6}. However, his methods were not powerful enough to yield the matrix analogue of Baker's conjecture, which recently surfaced as a conjecture posed by Kleinbock and Margulis in joint work with Victor Beresnevich, \cite[\67.1]{BKM}.
It is not obvious that the matrix analogue of Baker's conjecture follows from the Kleinbock--Margulis--Wang theorem (unlike Baker's original conjecture which followed directly from Sprind\v zuk's conjecture). The challenge here is to prove the strong nonplanarity of the manifold $\{X\oplus X^2\oplus \cdots\oplus X^n : X\in \MM_{m,m}\}$. Our paper bridges this gap, thus resolving the matrix version of Baker's conjecture in the affirmative.


\subsection{Diophantine approximation of matrices}
In what follows we recall the theorem of Kleinbock--Margulis--Wang regarding strongly extremal manifolds in the space of matrices. We begin by recalling the definition of strong extremality. Fix $\bothdim\in\N$, let $\MM = \MM_\bothdim$ denote the set of $\pdim\times \qdim$ matrices, and fix $\bfA\in\MM$. The \emph{exponent of irrationality} of $\bfA$ is defined as
\[
\omega(\bfA) = \limsup_{\substack{\qq\in\Z^\qdim\butnot\{\0\} \\ \pp\in\Z^\pdim}} \frac{-\log\|\bfA\qq - \pp\|}{\log\|\qq\|},
\]
and the \emph{exponent of multiplicative irrationality} is the number\Footnote{This definition agrees with the multiplicative approximation framework considered in \cite{KMW}, but not the one considered in \cite{KleinbockMargulis}; see comments after \cite[Proposition 4.1]{DFSU_GE1} for more details.}
\[
\omega_\mult(\bfA) = \limsup_{\substack{\qq\in\Z^\qdim\butnot\{\0\} \\ \pp\in\Z^\pdim}} \frac{-\log\prod_{i = 1}^\pdim |(\bfA\qq - \pp)_i|}{\log \prod_{j = 1}^\qdim |q_j|\vee 1}\cdot
\]
Note that $\omega_\mult(\bfA) \geq (\pdim/\qdim)\omega(\bfA)$. The matrix $\bfA$ is called \emph{very well approximable} if $\omega(\bfA) > \qdim/\pdim$, and \emph{very well multiplicatively approximable} if $\omega_\mult(\bfA) > 1$. Every very well approximable matrix is also very well multiplicatively approximable, and the sets of very well approximable and very well multiplicatively approximable matrices are both Lebesgue nullsets of full Hausdorff dimension. Finally, a measure $\mu$ on $\MM$ is called \emph{extremal} (resp. \emph{strongly extremal}) if it gives zero measure to the set of very well approximable (resp. very well multiplicatively approximable) vectors.

Now we recall the definition of strong nonplanarity, the hypothesis that Kleinbock--Margulis--Wang need in order to prove that a manifold is strongly extremal. To this end, we make the following definitions:

\begin{definition}[Cf. {\cite[p.3]{BKM}}, {\cite[p.7]{DFSU_GE1}}]
\label{definitionplucker}
The \emph{Pl\"ucker embedding} is the map which sends a matrix to the list of the determinants of its minors.
\end{definition}

\begin{definition}[Cf. {\cite[p.341]{KleinbockMargulis2}}, {\cite[(2.2)]{BKM}}]
\label{definitionstronglynonplanar}
A submanifold $M\subset \R^d$ is said to be \emph{nondegenerate} if whenever $\phi:U\to M$ is a coordinate chart, then for all $x\in U$ the components of the derivatives $\phi'(x),\phi''(x),\ldots$ span $\R^d$. A submanifold $M\subset \MM$ is said to be \emph{strongly nonplanar} if its image under the Pl\"ucker embedding is nondegenerate.\Footnote{This definition agrees with the one given in \cite{BKM} for all real-analytic manifolds, but the definitions disagree for some manifolds which are not real-analytic. The validity of Theorem \ref{theoremKMW} in the non-analytic category shows that our definition is the ``correct'' one.}
\end{definition}

\begin{remark}
A connected real-analytic manifold is nondegenerate if and only if it is not contained in any affine hyperplane, and strongly nonplanar if and only if it is not contained in the preimage of any affine hyperplane under the Pl\"ucker embedding.
\end{remark}

\begin{remark}
The reason for the adjective ``strong'' is that there is another notion of nonplanarity, \emph{weak nonplanarity}, which was introduced in \cite{BKM} and is sufficient to imply strong extremality. Since the matrix analogue of the Veronese curve turns out to be strongly nonplanar, in this paper there is no need to consider weak nonplanarity.
\end{remark}

We can now state the matrix analogue of the Kleinbock--Margulis theorem:

\begin{theorem}[\cite{KMW}, {\cite[Theorem 1.7]{DFSU_GE1}}]
\label{theoremKMW}
Every strongly nonplanar manifold is strongly extremal.
\end{theorem}

\begin{remark}
In \cite{KMW} this theorem was stated under the additional hypothesis that the manifold in question is real-analytic. This hypothesis is necessary if ``strongly nonplanar'' is defined as in \cite{BKM}, but if its definition is as given above, then the argument of \cite{KMW} works for non-analytic manifolds as well. Moreover, the validity of Theorem \ref{theoremKMW} for non-analytic manifolds also follows directly from \cite[Theorem 1.7]{DFSU_GE1}.
\end{remark}

\subsection{Main results}
Theorem \ref{theoremKMW} shows that the following theorem implies the matrix version of Baker's conjecture:

\begin{theorem}
\label{theorem1}
Fix $m,n\in\N$ and define the map $\psi_{m,n}:\MM_{m,m}\to \MM_{m,mn}$ via the formula\Footnote{In this paper $X\oplus Y$ denotes the direct columnwise sum of two matrices $X$ and $Y$, i.e. the matrix whose columns are those of $X$ followed by those of $Y$.}
\begin{equation}
\label{psidef1}
\psi_{m,n}(X) = X\oplus X^2\oplus\cdots\oplus X^n.
\end{equation}
Then the manifold $\psi_{m,n}(\MM_{m,m}) \subset \MM_{m,mn}$ is strongly nonplanar.
\end{theorem}

\begin{corollary}[Matrix version of Baker's conjecture]
For Lebesgue almost every $X\in\MM_{m,m}$, the matrix $\psi_{m,n}(X)$ defined by \eqref{psidef1} is not very well multiplicatively approximable.
\end{corollary}

The techniques we use to prove Theorem \ref{theorem1} also allow us to prove a more general result. Returning momentarily to the simultaneous approximation setting, a well-known generalization of the Veronese curve is the \emph{Veronese manifold} of dimension $\cardQ$ and order $n$, which is the image of $\R^\cardQ$ under the \emph{Veronese embedding}
\[
\psi_{n,\cardQ}(x_1,\ldots,x_\cardQ) \df \left(\prod_{k = 1}^\cardQ x_k^{n_k}\right)_{1 \leq \sum n_k \leq n}.
\]
The Veronese curve corresponds to the special case $\cardQ = 1$. It is readily verified that all Veronese manifolds are nondegenerate, so by the Kleinbock--Margulis theorem \cite{KleinbockMargulis2}, Veronese manifolds are strongly extremal. The following result implies that the same is true for the matrix versions of these manifolds:

\begin{theorem}
\label{theorem2}
Fix $m,n,\cardQ\in\N$, and for each $X_1,\ldots,X_\cardQ\in \MM_{m,m}$ let
\begin{equation}
\label{psidef2}
\psi_{m,n,\cardQ}(X_1,\ldots,X_\cardQ) \df \bigoplus_{1 \leq \sum n_k \leq n}\prod_{k = 1}^\cardQ X_k^{n_k}.
\end{equation}
Then $\psi_{m,n,\cardQ}(\MM_{m,m}^\cardQ)$ is strongly nonplanar (and thus strongly extremal).
\end{theorem}

One aspect of Theorem \ref{theorem2} which is somewhat undesirable is that due to the noncommutativity of matrices, the ordering of the product $\prod_{k = 1}^\cardQ X_k^{n_k}$ matters.\Footnote{On the other hand, the ordering of the direct sum $\bigoplus_{1\leq\Sigma n_k\leq n}$ does not matter, since applying a transformation rearranging the columns of a matrix does not affect the strong nonplanarity of a manifold; it is equivalent to applying a linear transformation to the image of the matrix under the Pl\"ucker embedding.} In Theorem \ref{theorem2} we understand the ordering to be the default one, that is, $\prod_{k = 1}^\cardQ X_k^{n_k} = X_1^{n_1}\cdots X_\cardQ^{n_\cardQ}$. But other orderings are possible, which would lead to different theorems. We now proceed to describe a general theorem that will cover all possible orderings.

Fix $\cardQ\in\N$, and let $M_\cardQ$ (resp. $M_\cardQ^*$) denote the set of all monomials in the commutative (resp. noncommutative) dummy variables $x_1,\ldots,x_\cardQ$. So e.g. $x_1 x_2$ and $x_2 x_1$ denote distinct elements of $M_\cardQ^*$, but they denote the same element of $M_\cardQ$. Let $\pi:M_\cardQ^*\to M_\cardQ$ be the natural map, i.e. the one that sends the element of $M_\cardQ^*$ denoted by $x_1 x_2$ to the element of $M_\cardQ$ denoted by $x_1 x_2$. Finally, given $\poly\in M_\cardQ^*$ and matrices $X_1,\ldots,X_\cardQ\in\MM_{m,m}$, let $\poly(X_1,\ldots,X_\cardQ)$ denote the matrix that results from substituting $x_1 = X_1,\ldots,x_\cardQ = X_\cardQ$ into $\poly$. So e.g. $(x_1 x_2 x_1)(X_1,X_2) = X_1 X_2 X_1$.

\begin{theorem}
\label{theorem3}
Fix $m,n,\cardQ\in\N$, and let $\bfP = (\poly_\ell)_{\ell = 1}^n$ be a finite sequence in $M_\cardQ^*$ such that $\pi(\poly_1),\ldots,\pi(\poly_n)$ are distinct. For each $X_1,\ldots,X_\cardQ\in \MM_{m,m}$, let
\begin{equation}
\label{psidef3}
\psi_{m,n,\cardQ}^{(\bfP)}(X_1,\ldots,X_\cardQ) = \bigoplus_{\ell = 1}^n \poly_\ell(X_1,\ldots,X_\cardQ).
\end{equation}
Then $\psi_{m,n,\cardQ}^{(\bfP)}(\MM_{m,m}^\cardQ)$ is strongly nonplanar (and thus strongly extremal).
\end{theorem}

Obviously, Theorem \ref{theorem1} is a special case of Theorem \ref{theorem2} and Theorem \ref{theorem2} is a special case of Theorem \ref{theorem3}. So in what follows, we will prove only Theorem \ref{theorem3}.

\begin{remark}
\label{remarknoninj}
In our proof, the requirement that $\pi(\poly_1),\ldots,\pi(\poly_n)$ be distinct is a necessary one. Indeed, if $\pi(\poly_1),\ldots,\pi(\poly_n)$ are not distinct, then we may choose $i\neq j$ such that $\pi(\poly_i) = \pi(\poly_j)$, and then for all $X_1,\ldots,X_\cardQ$ we have
\[
\det[\poly_i(X_1,\cdots,X_\cardQ)] = \det[\poly_j(X_1,\cdots,X_\cardQ)]
\]
which implies that $\psi_{m,n,\cardQ}^{(\bfP)}((\MM_{m,m})^\cardQ)$ is not strongly nonplanar (its image under the Pl\"ucker embedding is contained in a linear subspace of the form $x = y$, where $x$ and $y$ are coordinate functions). However, in this case the manifold $\psi_{m,n,\cardQ}^{(\bfP)}((\MM_{m,m})^\cardQ)$ might still be extremal. For example, if $m = n = \cardQ = 2$, $p_1 = x_1 x_2$, and $p_2 = x_2 x_1$, then it can be checked directly that $\psi_{m,n,\cardQ}^{(\bfP)}((\MM_{m,m})^\cardQ)$ is weakly nonplanar; thus by \cite[Corollary 2.4]{BKM}, it is strongly extremal. It is unclear whether all manifolds of the form $\psi_{m,n,\cardQ}^{(\bfP)}((\MM_{m,m})^\cardQ)$ are extremal.
\end{remark}


\section{Multisets}
To facilitate the proof of Theorem \ref{theorem3}, we introduce some notation for dealing with multisets. If $I$ is a set, then a \emph{multisubset} of $I$ is a map from $I$ to $\N \df \{0,1,\ldots\}$; a \emph{multiset} is a multisubset of some set. We use the notation $S\prec I$ to mean that $S$ is a multisubset of $I$, and we use the notation $\#(S;i)$ to denote the value of $S$ at a point $i\in I$. (Intuitively, $\#(S;i)$ is the ``number of copies of $i$ in $S$''.) The \emph{cardinality} of a multiset $S\prec I$ is $\#(S) \df \sum_{i\in I} \#(S;i)$. If $S\prec I$ and $(a_i)_{i\in I}$ is a collection of real numbers, then
\begin{align*}
\sum_{i\in S} a_i &\df \sum_{i\in I} {\#(S;i)}a_i,&
\prod_{i\in S} a_i &\df \prod_{i\in I} a_i^{\#(S;i)}.
\end{align*}
If $S\subset I$ is a subset, then $S$ can also be interpreted as a multiset via its characteristic function $\#(S;i) \df \#(S\cap\{i\})$, in which case the notations above agree with the usual definition of sums and products over $S$. Conversely, if $S\prec I$ is a multiset, then we define the \emph{support} of $S$ to be the unique set $|S|\subset I$ such that for all $i\in I$, we have $i\in |S|$ if and only if $\#(S;i) \geq 1$. Note that if $|S|$ is then reinterpreted as a multiset, it is possible that $S\neq |S|$; in fact, this happens whenever $\#(S;i) \geq 2$ for some $i\in I$.

If $S,T\prec I$, then we say that $S$ is a \emph{submultiset}\Footnote{To warn against confusing the two words, note that a multisubset of $I$ is a submultiset of $I$ if and only if it is a set.} of $T$, denoted $S\leq T$, if $\#(S;i) \leq \#(T;i)$ for all $i\in I$. If $S\prec I$ and $(T_i)_{i\in I}$ is a collection of (multi)subsets of a set $J$, then we define the multiset $\sum_{i\in S} T_i \prec J$ via the formula
\[
\#\left(\sum_{i\in S} T_i;j\right) \df \sum_{i\in S} \#(T_i;j) \all j\in J.
\]
Note that even if $S$ and $(T_i)_{i\in I}$ are all ordinary sets, the multiset $\sum_{i\in S} T_i$ may have points with multiplicity $\geq 2$.

\section{Proof of Theorem \ref{theorem3}}
For convenience, in the proof we let $\AA = (\MM_{m,m})^\cardQ$, $\BB = \MM_{m,mn}$, and $\CC = \R^c$, where $c = \sum_{k = 1}^m \binom mk \binom{mn}k$ is the number of minors that an $m\times mn$ matrix has. We let $\psi = \psi_{m,n,\cardQ}^{(\bfP)}:\AA\to\BB$ be defined as in \eqref{psidef3}, and we let $\phi:\BB\to\CC$ be the Pl\"ucker embedding (cf. Definition \ref{definitionplucker}). By definition, to show that the manifold $\psi(\AA) \subset \BB$ is strongly nonplanar, we need to show that the manifold $\phi\circ\psi(\AA)\subset \CC$ is nondegenerate. Since $\phi\circ\psi$ is real-analytic (in fact polynomial), it suffices to show that $\phi\circ\psi(\AA)$ is not contained in any affine hyperplane, and since $\phi\circ\psi(\0) = \0$, it suffices to consider linear hyperplanes. So we need to show that $\VV = \CC$, where $\VV$ denotes the linear span of $\phi\circ\psi(\AA)$.

We begin by defining appropriate bases of the vector spaces $\AA$, $\BB$, and $\CC$. Let $V = \{1,\ldots,m\}$, $Q = \{1,\ldots,\cardQ\}$, and $E = V\times V\times Q$, and for each $(i,t,q)\in E$ let $\ee_{i,t,q} \in \AA$ be the basis vector whose $q$th coordinate is the elementary matrix $\ee_{i,t}$ whose $(i,t)$th entry is 1. The letters $V$ and $E$ are chosen because later in the proof, we will think of $(V,E)$ as a directed multigraph, where the edge $(i,t,q)\in E$ is considered to have initial vertex $i$ and terminal vertex $t$. Let $L = \{1,\ldots,n\}$ and $W = \{1,\ldots,mn\}$, and for each $t\in V$ and $\ell\in L$ let $\lb t,\ell\rb = (\ell - 1) m + t \in W$. For each $(i,j) \in V\times W$ let $\ee_{i,j}\in\BB$ be the elementary matrix whose $(i,j)$th entry is $1$, i.e.
\[
\ee_{i,\lb t,\ell\rb} \df \0\oplus\cdots\oplus\underbrace{\ee_{i,t}}_{\text{$\ell$th}}\oplus\cdots\oplus\0.
\]
Finally, let
\[
D = \{(I,J) : I\subset V,\; J\subset W, \; \#(I) = \#(J) \geq 1\},
\]
and let $\theta:D\leftrightarrow\{1,\ldots,c\}$ be the bijection that describes the order in which the Pl\"ucker embedding $\phi$ lists the determinants of the minors. For each $(I,J)\in D$, let $\ee_{I,J}\in \R^c$ be the basis vector whose $\theta(I,J)$th coordinate is $1$.

Given $A\in \AA$, $B\in \BB$, and $C\in \CC$, we will denote by $A_{i,t,q}$, $B_{i,j}$, and $C_{I,J}$ the coordinates of $A$, $B$, and $C$ with respect to the bases defined above, so that
\begin{align*}
A &= \sum_{(i,t,q)\in E} A_{i,t,q} \ee_{i,t,q},&
B &= \sum_{(i,j)\in V\times W} B_{i,j} \ee_{i,j},&
C &= \sum_{(I,J)\in D} C_{I,J} \ee_{I,J}.
\end{align*}
Now let $B = \psi(A)$ and $C = \phi(B)$. By the definition of the Pl\"ucker embedding $\phi$, $C_{I,J}$ is the determinant of the matrix $(B_{i,j})_{i\in I, j\in J}$, i.e.\Footnote{If $I,J$ are ordered sets, the \emph{sign} of a bijection $\sigma:I\leftrightarrow J$ is defined as $\sgn(\sigma)\df (-1)^{\#\{(i_1,i_2)\in I^2:i_1 < i_2,\;\sigma(i_1) > \sigma(i_2)\}}$.}
\begin{equation}
\label{CIJ}
C_{I,J} = \sum_{\sigma:I\leftrightarrow J} \sgn(\sigma) \prod_{i\in I} B_{i,\sigma(i)}.
\end{equation}
Similarly, by the definition of $\psi$,
\begin{equation}
\label{Bitl}
B_{i,\lb t,\ell\rb} = [\poly_\ell(X_1,\ldots,X_\cardQ)]_{i,t},
\end{equation}
where $A = (X_1,\ldots,X_\cardQ)$.

For each $\ell\in L$ let $\|\ell\|\in\N$ and $f_\ell(1),\ldots,f_\ell(\|\ell\|)\in Q$ be selected so that
\[
\poly_\ell(X_1,\ldots,X_\cardQ) = X_{f_\ell(1)}\cdots X_{f_\ell(\|\ell\|)}.
\]
Note that these choices are unique because $\poly_\ell$ is a monomial in noncommutative variables. We will later need the following fact about the commutative version of $\poly_\ell$: since it can be written in multiset notation as
\[
\pi(\poly_\ell) = \prod_{q\in \sum_{k = 1}^{\|\ell\|} \{f_\ell(k)\}} x_q,
\]
it follows that $\pi(\poly_\ell)$ depends only on $\sum_{k = 1}^{\|\ell\|} \{f_\ell(k)\}$.

The next step is to write the coordinates $B_{i,\lb t,\ell\rb} = [X_{f_\ell(1)}\cdots X_{f_\ell(\|\ell\|)}]_{i,t}$ ($i,t\in V$, $\ell\in L$) of $B$ in terms of the coordinates $A_{i,t,q} = [X_q]_{i,t}$ ($(i,t,q)\in E$) of $A$ using the definition of matrix multiplication. To facilitate this, we define a \emph{path} in $(V,E)$ to be a diagram of the form
\begin{equation}
\label{pathdiagram}
v_0 \tendsto{f_\ell(1)} v_1 \tendsto{f_\ell(2)} \cdots \tendsto{f_\ell(\|\ell\|)} v_{\|\ell\|},
\end{equation}
where $\ell\in L$ and $v_0,\ldots,v_{\|\ell\|}\in V$. Let $P$ denote the path represented by the diagram \eqref{pathdiagram}. Then the vertices $i(P) \df v_0$ and $t(P) \df v_{\|\ell\|}$ are called the \emph{initial} and \emph{terminal vertices} of $P$, and $\ell(P) \df \ell$ is called the \emph{label} of $P$. The triples $e_k(P) \df (v_{k - 1},v_k,f_\ell(k)) \in E$ ($k = 1,\ldots,\|\ell\|$) are called the \emph{edges} of $P$, and the multiset $F(P) \df \sum_{k = 1}^{\|\ell\|} \{e_k(P)\}$ is called the \emph{edge multiset} of $P$.

Using this notation, we can rewrite \eqref{Bitl} as
\[
B_{i,\lb t,\ell\rb} = [\poly_\ell(X_1,\ldots,X_\cardQ)]_{i,t} = \sum_P \prod_{k = 1}^{\|\ell\|} [X_{f_\ell(k)}]_{v_{k - 1},v_k} = \sum_P \prod_{e\in F(P)} A_e,
\]
where the sums are taken over all paths $P$ such that $i(P) = i$, $t(P) = t$, and $\ell(P) = \ell$. Plugging into \eqref{CIJ} and using the notation
\begin{align*}
j(P) &\df \lb t(P),\ell(P)\rb \in W,&
I(\PP) &\df \sum_{P\in\PP} \{i(P)\} \prec V,\\
J(\PP) &\df \sum_{P\in\PP} \{j(P)\} \prec W,&
F(\PP) &\df \sum_{P\in\PP} F(P) \prec E,\\
d(\PP) &\df (I(\PP),J(\PP)).
\end{align*}
gives
\begin{align*}
C_{I,J}
&= \sum_{\sigma:I\leftrightarrow J} \sgn(\sigma) \prod_{i\in I} \sum_{\substack{P \\ i(P) = i \\ j(P) = \sigma(i)}} \prod_{e\in F(P)} A_e\\
&= \sum_{\sigma:I\leftrightarrow J} \sgn(\sigma) \sum_{\substack{\PP \\ I(\PP) = I \\ j(P) = \sigma(i(P)) \all P\in |\PP|}} \prod_{P\in\PP} \prod_{e\in F(P)} A_e\\
&= \sum_{\substack{\PP \\ d(\PP) = (I,J)}} (\pm 1) \prod_{e\in F(\PP)} A_e,
\end{align*}
where the last two sums are taken over all multisubsets $\PP$ of the collection of paths. Note that every multiset $\PP$ included in the summation, i.e. for which $I(\PP) = I$ and $J(\PP) = J$, is in fact a set (by which we mean that it satisfies $\PP = |\PP|$). Also note that the inclusion $d(\PP) \in D$ is valid for all $\PP$ that appear in the summation, but not for every collection of paths $\PP$.

\begin{lemma}
\label{lemmapartialorder}
There exists a partial order $<$ on $D$ such that for all $d\in D$, there exists a collection of paths $\PP$ such that $d(\PP) = d$ and such that if $\PP'$ is another collection of paths such that $d(\PP')\in D$ and $F(\PP') = F(\PP)$, then $d(\PP') < d$.
\end{lemma}
\begin{subproof}
For each $d = (I,J) \in D$ let
\[
f(d) = \sum_{\substack{i\in I \\ J[i] \neq \smallemptyset}} \min_{\substack{\ell\in J[i]}} \|\ell\|,
\]
where $J[i] \df \{\ell\in L : \lb i,\ell\rb\in J\}$ for each $i\in I$. The partial order on $D$ will be defined as follows: we write $d_1 < d_2$ if and only if $f(d_1) < f(d_2)$. Now fix $d = (I,J)\in D$, let $S = \{i\in I : J[i] \neq \emptyset\}$, and for each $i\in S$ let $\ell_i \in J[i]$ be chosen to minimize $\|\ell_i\|$, so that $f(d) = \sum_{i\in S} \|\ell_i\|$. Let $\sigma':S\to J$ be defined by the equation $\sigma'(i) = \lb i,\ell_i\rb$, and let $\sigma:I\leftrightarrow J$ be an arbitrary bijective extension of $\sigma'$. For each $i\in I$ write $\sigma(i) = \lb t_i,\ell_i\rb$ (note that this agrees with the previous definition of $\ell_i$ if $i\in S$) and let $\path_i$ be the path represented by the diagram
\[
i \tendsto{f_{\ell_i}(1)} i\tendsto{f_{\ell_i}(2)} \cdots\tendsto{f_{\ell_i}(\|\ell_i\| - 1)}i \tendsto{f_{\ell_i}(\|\ell_i\|)} t_i,
\]
so that $i(\path_i) = i$ and $j(\path_i) = \sigma(i)$. Then, letting $\PP = \{\path_i : i\in I\}$, we have $I(\PP) = I$ and $J(\PP) = \sigma(I) = J$, i.e. $d(\PP) = d$.

Now let
\[
F = F(\PP) = \sum_{i\in I} \left[\sum_{k = 1}^{\|\ell_i\| - 1} \{(i,i,f_{\ell_i}(k))\} + \{(i,t_i,f_{\ell_i}(\|\ell_i\|))\}\right],
\]
and let $\PP'\neq\PP$ be a collection of paths such that $F(\PP') = F$ and $d(\PP')\in D$. Now the only edges $(i,t,q)\in |F|$ such that $i\neq t$ are those of the form $e_i \df (i,t_i,f_{\ell_i}(\|\ell_i\|))$ ($i\in I\butnot S$), which form a directed graph that contains no paths of length greater than one. This implies that if $\path_i'$ denotes the unique element of $\PP'$ such that $e_i\in |F(\path_i')|$, then:
\begin{itemize}
\item for all $i\in I\butnot S$, the diagram representing $\path_i'$ is of the form
\begin{equation}
\label{iti}
i\to\cdots\to i\to t_i\to\cdots\to t_i,
\end{equation}
and in particular $i(\path_i') = i$ and $t(\path_i') = t_i$;
\item for all $\path'\in \QQ \df \PP'\butnot\{\path_i' : i\in I\butnot S\}$, the diagram representing $\path'$ is of the form
\begin{equation}
\label{iii}
i\to\cdots\to i,
\end{equation}
where $i = i(\path_i) = t(\path_i)$.
\end{itemize}
Now since $d(\PP')\in D$, $I(\PP')$ is a set. Thus the paths $P_i'$ ($i\in I\butnot S$) are distinct, and
\begin{equation}
\label{Piprimepushout}
i(\path') \notin I\butnot S \text{ for all } \path'\in\QQ.
\end{equation}
It follows from \eqref{Piprimepushout} that $t(\path')\notin I\butnot S$ for all $\path'\in\PP'$, so we have $J'[i] = \emptyset$ for all $i\in I\butnot S$, where we write $d(\PP') = (I',J')$. So
\[
\{i'\in I' : J'[i'] \neq \emptyset\} = \{i(\path') : \path'\in |\QQ|\} \leq \sum_{\path' \in \QQ} \{i(\path')\}
\]
and thus
\[
f(d(\PP')) \leq \sum_{\path'\in\QQ} \min_{\ell\in J'[i(\path')]} \|\ell\| \leq \sum_{\path'\in\QQ} \|\ell(\path')\| = \sum_{\path'\in\QQ} \#(F(\path')) = \#(F) - \sum_{i\in I\butnot S} \#(F(\path_i')).
\]
On the other hand, combining \eqref{Piprimepushout} with the fact that all edges of $P_i$ have initial vertex $i$ gives that $F(\path_i) \leq F(\path_i')$ for all $i\in I\butnot S$. So
\[
f(d(\PP')) \leq \#(F) - \sum_{i\in I\butnot S} \#(F(\path_i')) \leq \#(F) - \sum_{i\in I\butnot S} \#(F(\path_i)) = f(d(\PP)).
\]
If strict inequality holds, we are done, so suppose equality holds. Then $F(\path_i) = F(\path_i')$ for all $i\in I\butnot S$. Since $\sum_{P'\in\QQ} F(P') = \sum_{i\in S} \sum_{k = 1}^{\|\ell_i\|} \{(i,i,f_{\ell_i}(k))\}$, the injectivity of the map $\PP'\ni \path'\mapsto i(\path')$ implies that we can write $\QQ = \{\path_i' : i\in S\}$ for some paths $\path_i'$ ($i\in S$) satisfying $F(\path_i) = F(\path_i')$. So $\PP' = \{\path_i' : i\in I\}$, and $F(\path_i) = F(\path_i')$ for all $i\in I$. We claim that $\path_i = \path_i'$ for all $i\in I$. Indeed, fix $i\in I$. Due to the diagram forms \eqref{iti} and \eqref{iii}, the paths $\path_i$ and $\path_i'$ agree except possibly for their labels $\ell = \ell(\path_i)$ and $\ell' = \ell(\path_i')$. But since $F(\path_i) = F(\path_i')$, we have $\sum_{k = 1}^\ell \{f_\ell(k)\} = \sum_{k = 1}^{\ell'} \{f_{\ell'}(k)\}$. But since the polynomials $\pi(\poly_1),\ldots,\pi(\poly_n)$ are assumed to be distinct, this implies that $\ell = \ell'$ and thus that $\path_i = \path_i'$. Since $i$ was arbitrary, we get $\PP' = \PP$, a contradiction.
\end{subproof}

To finish the proof, we will show that $\ee_d\in\VV$ for all $d\in D$ by using strong induction on the partial order $<$. (Recall that $\VV$ denotes the linear span of $\phi\circ\psi(\AA)$, for which we must show that $\VV = \CC$.) Fix $d\in D$ and suppose that for all $d' < d$, we have $\ee_{d'}\in \VV$. Let the collection of paths $\PP$ be as in Lemma \ref{lemmapartialorder}, and let $F = F(\PP)$. Let $\vv_F\in\CC$ be the coefficient of $\prod_{e\in F} A_e$ that appears in $\phi\circ \psi(A)$, thought of as a polynomial function of $A$. Then we have
\[
\VV\ni \vv_F = \sum_{\substack{\PP' \\ F(\PP') = F}} (\pm 1) \ee_{d(\PP')} = (\pm 1) \ee_d + \sum_{d' < d} \sum_{\substack{\PP' \\ d(\PP') = d' \\ F(\PP') = F}} (\pm 1) \ee_{d'}.
\]
Since by induction $\ee_{d'} \in \VV$ for all $d' < d$, solving for $\ee_d$ gives $\ee_d\in \VV$, completing the proof.

\section{Concluding remarks and further questions}

\subsection{An algorithm for computing weak nonplanarity}
\label{subsectionalgorithm}

In \cite[\67.1]{BKM}, the matrix version of Baker's conjecture is motivated by the fact that it is a ``special case'' of the question of ``how one can in general show that a given [manifold] is weakly nonplanar''. This may appear to be a broadly philosophical question to which there perhaps is no ultimate precise answer. However, it does motivate the following mathematical question that can be regarded as a certain way of making the philosophical question precise: is there an algorithm for determining whether a given manifold is weakly nonplanar, based on input parameters which describe the manifold uniquely in some way (e.g. via a parameterization)? In the next paragraph, we prove the existence of such an algorithm under the assumption that the manifold is parameterized by an algebraic mapping. 
That said, one should be warned that this result is not very useful for practical purposes. In particular, it cannot be used to prove the matrix version of Baker's conjecture, since that conjecture has an infinite list of submanifolds that need to be proven nonplanar, and it is impossible to run the algorithm on all those manifolds simultaneously. Still, the algorithmically checkable nature of weak nonplanarity it is an interesting theoretical result in its own right, complementing the main results of this note.

Recall (cf. \cite[(2.5)]{BKM}) that a connected real-analytic manifold $M\subset\MM = \MM_{m,n}$ is said to be \emph{weakly nonplanar} if it is not contained in any set of the form
\[
\HH_{A,B} = \{Y\in\MM : \det(AY + B) = 0\},
\]
where $A\in\MM_{n,m}$ and $B\in\MM_{n,n}$ are not both zero. It is shown in \cite[Corollary 2.4]{BKM} that every weakly nonplanar manifold is strongly extremal. Now suppose that the manifold $M$ is given in the form $M = \ff(\R^k)$, where $\ff:\R^k\to \MM_{m,n}$ is a polynomial map. Then $M$ is weakly nonplanar if and only if there exist $A,B$, not both zero, such that for all $\xx\in\R^k$, $\det(A \ff(\xx) + B) = 0$. This condition is expressible in first-order quantifier logic over the real numbers. Thus there is an algorithm \cite[Algorithm 11.14]{BPR} to determine whether this condition is true or not.

We remark that if $A,B$ are allowed to be complex matrices rather than real ones, then there is an algorithm for determining the truth of the condition which is somewhat simpler than the one described in \cite[Algorithm 11.14]{BPR}. In practice this is probably the better algorithm to use, since the case where $M$ is weakly nonplanar even though it is contained in a manifold of the form $\HH_{A,B}$ with $A,B$ complex matrices is rare.\Footnote{It happens for example when \begin{align*} M &= \left\{\left[\begin{array}{cc} x & y \\ -y & x\end{array}\right] : x,y\in\R\right\},& A &= \left[\begin{array}{cc} 1 & \sqrt{-1} \\ 0 & 0 \end{array}\right],& B &= \left[\begin{array}{cc} 0 & 0 \\ 1 & \sqrt{-1} \end{array}\right]. \end{align*}} For this algorithm, first write $\det(A \ff(\xx) + B) = \gg(A,B)\cdot \psi\circ\ff(\xx)$ for some polynomial function $\gg$, where $\psi$ is the Pl\"ucker embedding. Then compute the linear span of $\psi\circ\ff(\R^k)$; this is possible since $\psi\circ\ff$ is a polynomial. Write the linear span in the form $\lb \zz_1,\ldots,\zz_n\rb$. Then the question is whether there exist complex matrices $A,B$, not both zero, such that
\[
\gg(A,B)\cdot \zz_i = 0 \all i = 1,\ldots,n.
\]
This is a finite number of polynomial equations in the entries of the matrices $A$ and $B$. Hilbert's Nullstellensatz implies that the existence of such matrices is equivalent to the condition that there is some entry $x$ of $A$ or $B$ such that no power of $x$ is in the ideal $I$ generated by the polynomials $\gg(A,B)\cdot\zz_i$ ($i = 1,\ldots,n$). Given any entry $x$, it can be determined whether this condition holds by computing a Gr\"obner basis for $I$ with respect to some monomial ordering with respect to which the monomials $x,x^2,\ldots$ are ordered lower than all other monomials, and then checking whether any power of $x$ is in the Gr\"obner basis.

\subsection{Algorithms for computing extremality and strong extremality within the class of rational submanifolds}
Subsequent to the question asked in \cite[\67.1]{BKM}, Aka, Breuillard, Rosenzweig, and de Saxc\'e found a sufficient condition for extremality which is weaker than weak nonplanarity \cite[Theorem 3.1]{ABRdS}. Their result is optimal in the sense that it completely characterizes extremality within the class of rational submanifolds (i.e. those defined by polynomial equations with rational coefficients), indicating that extremal manifolds which do not satisfy their condition must be extremal for ``Diophantine'' reasons rather than ``geometric'' ones. The above argument goes through with minor modifications to show that there is an algorithm for determining whether their condition holds for any given manifold. In particular, there is an algorithm for determining whether any rational submanifold of $\MM$ is extremal or not.

A forthcoming joint work of the authors with Breuillard and de Saxc\'e will include the analogous version of their result for strong extremality. The above argument would also apply to this setting, producing an algorithm for determining whether any rational submanifold of $\MM$ is strongly extremal or not.

\subsection{Open question}
We end by posing the following question, which would extend the main result of our paper, and which we are unable to resolve at present.
\begin{question}[Cf. Remark \ref{remarknoninj}]
Fix $m,n,\cardQ\in\N$ and let $(\poly_\ell)_{\ell = 1}^n$ be a finite sequence in $M_r^*$, with $\pi(\poly_1),\ldots,\pi(\poly_n)$ not necessarily distinct. Is it true that the manifold $\psi_{m,n,\cardQ}^{(\bfP)}((\MM_{m,m})^\cardQ)$ must be
\begin{itemize}
\item[(i)] weakly nonplanar?
\item[(ii)] strongly extremal?
\item[(iii)] extremal?
\end{itemize}
(The distinctness of $\pi(\poly_1),\ldots,\pi(\poly_n)$ is a necessary and sufficient condition for the strong nonplanarity of $\psi_{m,n,\cardQ}^{(\bfP)}((\MM_{m,m})^\cardQ)$, as shown by combining Theorem \ref{theorem3} with Remark \ref{remarknoninj}.)
\end{question}

{\bf Acknowledgements.} The authors thank Nicolas de Saxc\'e for helpful comments on an earlier version of this paper. The authors also thank Chris Miller for pointing us to the reference \cite{BPR}.
The first-named author was supported in part by a 2016-2017 Faculty Research Grant from the University of Wisconsin-La Crosse. The research of the second-named author was supported in part by the EPSRC Programme Grant EP/J018260/1.

\bibliographystyle{amsplain}

\bibliography{bibliography}

\end{document}